\documentclass[a4paper,11pt,leqno]{article}
\usepackage[T1]{fontenc}
\usepackage{amssymb,latexsym}
\usepackage{latexsym}
\usepackage{amsfonts}
\usepackage{eucal}
\usepackage{amsmath}

\input{tcilatex}

\begin{document}

\title{Refinements of Hermite-Hadamard inequality on simplices}
\date{}
\author{\textbf{Flavia-Corina Mitroi and C\u{a}t\u{a}lin Irinel Spiridon}}
\maketitle

{\scriptsize Some refinements of the Hermite-Hadamard inequality are
obtained in the case of continuous convex functions defined on simplices. 
\newline
}

{\scriptsize \textit{AMS 2010 Subject Classification:} 26A51. \newline
}

{\scriptsize \textit{Key words:} Hermite-Hadamard inequality, convex
function, simplex.\newline
}

\begin{center}
{\normalsize \setcounter{section}{1}\textbf{1.INTRODUCTION} }
\end{center}

{\normalsize The aim of this paper is to provide a refinement of the
Hermite-Hadamard inequality on simplices. }

{\normalsize Suppose that $K$ is a metrizable compact convex subset of a
locally convex Hausdorff space $E.$ Given a Borel probability measure $\mu\ $%
on $K,$ one can prove the existence of a unique point $b_{\mu}\in K$ (called
the \emph{barycenter }of $\mu)$ such that 
\begin{equation*}
x^{\prime}(b_{\mu})=\int_{K}\,x^{\prime}(x)\,\mathrm{d}\mu(x)
\end{equation*}
for all continuous linear functionals $x^{\prime}$ on $E$. The main feature
of barycenter is the inequality 
\begin{equation*}
f(b_{\mu})\leq\int_{K\,}\,f(x)\,\mathrm{d}\mu(x),
\end{equation*}
valid for every continuous convex function $f:K\rightarrow\mathbb{R}$. It
was noted by several authors that this inequality is actually equivalent to
the Jensen inequality. }

{\normalsize The following theorem due to G. Choquet complements this
inequality and relates the geometry of $K$ to a given mass distribution.%
\newline
}

{\normalsize THEOREM 1. \label{choquet}$($The general form of
Hermite-Hadamard inequality$)$. Let $\mu$ be a Borel probability measure on
a metrizable compact convex subset $K$\ of a locally convex Hausdorff
space.\ Then there exists a Borel probability measure $\nu$\ on $K$ which
has the same barycenter as $\mu,$\ is zero outside $\limfunc{Ext}$\thinspace$%
K,$ and verifies the double inequality%
\begin{equation}
f(b_{\mu})\leq\int_{K\,}\,f(x)\,\mathrm{d}\mu(x)\leq\int_{\limfunc{Ext}%
K\,}\,f(x)\,\mathrm{d}\nu(x)  \label{HH}
\end{equation}
for all continuous convex functions $f:K\rightarrow\mathbb{R}.$ }

{\normalsize Here $\limfunc{Ext}$\thinspace$K$ denotes the set of all
extreme points of $K$. }

{\normalsize The details can be found in \cite{cpn06}, pp. 192-194. See also 
\cite{cpn02}. }

{\normalsize In the particular case of simplices one can take advantage of
the barycentric coordinates. }

{\normalsize Suppose that $\Delta\subset\mathbb{R}^{n}$ is an $n$%
-dimensional simplex of vertices $P_{1},...,P_{n+1}.$ In this case $\limfunc{%
Ext}\Delta =\{P_{1},...,P_{n+1}\}$ and each $x\in\Delta$ can be represented
uniquely as a convex combination of vertices,%
\begin{equation}
\sum_{k=1}^{n+1}\lambda_{k}\left( x\right) P_{k}\,=x,  \label{bar1}
\end{equation}
where the coefficients $\lambda_{k}\left( x\right) $ are nonnegative numbers
(depending on $x)$ and 
\begin{equation}
\sum_{k=1}^{n+1}\lambda_{k}\left( x\right) =1.  \label{bar2}
\end{equation}
Each function $\lambda_{k}:x\rightarrow\lambda_{k}\left( x\right) $ is an
affine function on $\Delta.$ This can be easily seen by considering the
linear system consisting of the equations (\ref{bar1}) and (\ref{bar2}). }

{\normalsize The coefficients $\lambda_{k}\left( x\right) $ can be computed
in terms of Lebesgue volumes (see \cite{bes08}, \cite{neu89}). We denote by $%
\Delta_{j}(x)$ the subsimplex obtained when the vertex $P_{j}$ is replaced
by $x\in\Delta.$ Then one can prove that 
\begin{equation}
\lambda_{k}\left( x\right) =\frac{\limfunc{Vol}\left( \Delta _{k}(x)\right) 
}{\limfunc{Vol}\left( \Delta\right) }  \label{lambda}
\end{equation}
for all $k=1,...,n+1$ (the geometric interpretation is very intuitive)$.$
Here $\limfunc{Vol}\left( \Delta\right) =\int_{\Delta}\mathrm{d}x.$ }

{\normalsize In the particular case where $\mathrm{d}\mu\left( x\right) =%
\mathrm{d}x/\limfunc{Vol}\left( \Delta\right) $ we have $\lambda_{k}\left(
b_{\mathrm{d}x/\limfunc{Vol}\left( \Delta\right) }\right) =\frac {1}{n+1}$
for every $k=1,...,n+1.$ }

{\normalsize The above discussion leads to the following form of Theorem \ref%
{choquet} in the case of simplices:\newline
}

{\normalsize COROLLARY 1. \label{cor1}Let $\Delta\subset\mathbb{R}^{n}$ be
an $n$-dimensional simplex of vertices $P_{1},...,P_{n+1}$ and $\mu$ be a
Borel probability measure on $\Delta$ with barycenter $b_{\mu}$. Then for
every continuous convex function $f:\Delta\rightarrow\mathbb{R}$,%
\begin{equation*}
f\left( b_{\mu}\right) \leq\int_{\Delta}\,f(x)\,\mathrm{d}\mu\left( x\right)
\leq\sum_{k=1}^{n+1}\lambda_{k}(b_{\mu})f\left( P_{k}\right) .
\end{equation*}
}

{\normalsize \textit{Proof}. In fact, 
\begin{align*}
\int_{\Delta}\,f(x)\,\mathrm{d}\mu\left( x\right) & =\int_{\Delta }\,f\left(
\sum_{k=1}^{n+1}\lambda_{k}\left( x\right) P_{k}\,\right) \mathrm{d}%
\mu\left( x\right) \leq\int_{\Delta}\sum_{k=1}^{n+1}\lambda _{k}\left(
x\right) \,f(P_{k})\,\mathrm{d}\mu\left( x\right) \\
& =\sum_{k=1}^{n+1}\,f(P_{k})\int_{\Delta}\lambda_{k}\left( x\right) \,\,%
\mathrm{d}\mu\left( x\right) =\sum_{k=1}^{n+1}\lambda_{k}(b_{\mu })f\left(
P_{k}\right) .
\end{align*}
On the other hand%
\begin{equation*}
\sum_{k=1}^{n+1}\lambda_{k}(b_{\mu})f\left( P_{k}\right) =\int _{\limfunc{Ext%
}\Delta}f\mathrm{d}\left( \sum_{k=1}^{n+1}\lambda
_{k}(b_{\mu})\delta_{P_{k}}\right)
\end{equation*}
and $\sum_{k=1}^{n+1}\lambda_{k}(b_{\mu})\delta_{P_{k}}$ is the only Borel
probability measure $\nu$ concentrated at the vertices of $\Delta$ which
verifies the inequality%
\begin{equation*}
\int_{\Delta\,}\,f(x)\,\mathrm{d}\mu(x)\leq\int_{\limfunc{Ext}\Delta
\,}\,f(x)\,\mathrm{d}\nu(x)
\end{equation*}
for every continuous convex functions $f:\Delta\rightarrow\mathbb{R}$.
Indeed $\nu$ must be of the form $\nu=\sum_{k=1}^{n+1}\alpha_{k}%
\delta_{P_{k}},$ with $b_{\nu}=\sum_{k=1}^{n+1}\alpha_{k}P_{k}=b_{\mu}$ and
the uniqueness of barycentric coordinates yields the equalities%
\begin{equation*}
\alpha_{k}=\lambda_{k}(b_{\mu})\text{ for }k=1,...,n+1.
\end{equation*}%
\newline
}

{\normalsize It is worth noticing that the Hermite--Hadamard inequality is
not just a consequence of convexity, it actually characterizes it. See \cite%
{tri08}. }

{\normalsize The aim of the present paper is to improve the result of
Corollary \ref{cor1}, by providing better bounds for the arithmetic mean of
convex functions defined on simplices. }

\begin{center}
{\normalsize \setcounter{section}{2}\textbf{2.MAIN RESULTS} }
\end{center}

{\normalsize We start by extending the following well known inequality
concerning the continuous convex functions defined on intervals: 
\begin{equation*}
\frac{1}{b-a}\int_{a}^{b}\,f(x)\,\mathrm{d}x\leq\frac{1}{2}\left( \frac{%
f\left( a\right) +f\left( b\right) }{2}+f\left( \frac{a+b}{2}\right) \right)
.
\end{equation*}
See \cite{cpn06}, p. 52.\newline
}

{\normalsize THEOREM 2. \label{thm2}Let $\Delta\subset\mathbb{R}^{n}$ be an $%
n$-dimensional simplex of vertices $P_{1},...,P_{n+1},$ endowed with the
normalized Lebesgue measure $\mathrm{d}x/\limfunc{Vol}\left( \Delta\right) .$
Then for every continuous convex function $f:\Delta\rightarrow\mathbb{R}$
and every point $P\in\Delta$ we have 
\begin{equation}
\frac{1}{\limfunc{Vol}\left( \Delta\right) }\int_{\Delta }\,f(x)\,\mathrm{d}%
x\leq\frac{1}{n+1}\left( \sum_{k=1}^{n+1}\left( 1-\lambda_{k}\left( P\right)
\right) f\left( P_{k}\,\right) +f\left( P\right) \right)  \label{rhh1}
\end{equation}
In particular, 
\begin{equation}
\frac{1}{\limfunc{Vol}\left( \Delta\right) }\int_{\Delta }\,f(x)\,\mathrm{d}%
x\leq\frac{1}{n+1}\left( \frac{n}{n+1}\sum_{k=1}^{n+1}f\left( P_{k}\right)
+f\left( b_{\mathrm{d}x/\limfunc{Vol}\left( \Delta\right) }\right) \right) .
\label{15}
\end{equation}
}

{\normalsize \textit{Proof}. Consider the barycentric representation $%
P=\sum_{k}\lambda_{k}\left( P\right) P_{k}\in\Delta$ . According to
Corollary \ref{cor1}, 
\begin{equation*}
\frac{1}{\lambda_{i}\left( P\right) \limfunc{Vol}\left( \Delta\right) }%
\int_{\Delta_{i}\left( P\right) }\,f(x)\,\mathrm{d}x\leq\frac{1}{n+1}\left(
\sum_{k\neq i}f\left( P_{k}\,\right) +f\left( P\right) \right) ,
\end{equation*}
where $\Delta_{i}\left( P\right) $ denotes the simplex obtained from $\Delta$
by replacing the vertex $P_{i}$ by $P$. Notice that (\ref{lambda}) yields $%
\lambda_{i}\left( P\right) \limfunc{Vol}\left( \Delta\right) =\limfunc{Vol}%
\left( \Delta_{i}\left( P\right) \right) .$ Multiplying both sides by $%
\lambda_{i}\left( P\right) $ and summing up over $i$ we obtain the
inequality (\ref{rhh1}). }

{\normalsize In the particular case when $P=b_{\mathrm{d}x/\limfunc{Vol}%
\left( \Delta\right) }$ all coefficients $\lambda_{k}\left( P\right) $ equal 
$1/(n+1)$. \newline
}

{\normalsize REMARK 1. Using \cite[Theorem 1]{bes08} instead of Corollary 1
of the present paper, one may improve the statement of Theorem \ref{thm3} by
the cancellation of the continuity condition. Such statement is proved
independently, by a different approach, in \cite{was11}.\newline
}

{\normalsize An extension of Theorem \ref{thm2} is as follows:\newline
}

{\normalsize THEOREM 3. \label{thm3}Let $\Delta\subset\mathbb{R}^{n}$ be an $%
n$-dimensional simplex of vertices $P_{1},...,P_{n+1}$ endowed with the
Lebesgue measure and let $\Delta^{\prime}\subseteq\Delta$ a subsimplex of
vertices $P_{1}^{\prime },...,P_{n+1}^{\prime}$ which has the same
barycenter as $\Delta$. Then for every continuous convex function $%
f:\Delta\rightarrow\mathbb{R},$ and every index $j\in\{1,...,n+1\},$ we have
the estimates%
\begin{align}
\frac{1}{\limfunc{Vol}\left( \Delta\right) }\int_{\Delta }\,f(x)\,\mathrm{d}%
x & \leq\frac{1}{n+1}\left( \sum_{k\neq j}\sum _{i}\lambda_{i}\left(
P_{k}^{\prime}\right) f\left( P_{i}\,\right) +f\left( P_{j}^{\prime}\right)
\right)  \label{R} \\
& \leq\frac{1}{n+1}\sum_{i}f\left( P_{i}\right) \   \notag
\end{align}
and%
\begin{align}
\frac{1}{\limfunc{Vol}\left( \Delta\right) }\int_{\Delta }\,f(x)\,\mathrm{d}%
x & \geq\sum_{i}\lambda_{i}\left( P_{j}^{\prime}\right) f\left( \frac{1}{n+1}%
\left( \sum_{k\neq i}P_{k}\,+P_{j}^{\prime}\right) \right)  \label{L} \\
& \geq f\left( \frac{1}{n+1}\sum_{i}P_{i}\right)  \notag
\end{align}
}

{\normalsize \textit{Proof}. We will prove here only the inequalities (\ref%
{R}), the proof of (\ref{L}) being similar. }

{\normalsize By Corollary \ref{cor1}, for each index $i,$ 
\begin{equation*}
\frac{1}{\lambda_{i}\left( P_{j}^{\prime}\right) \limfunc{Vol}\left(
\Delta\right) }\int_{\Delta_{i}\left( P_{j}^{\prime}\right) }f\left(
x\right) \mathrm{d}x\leq\frac{1}{n+1}\left( \sum_{k\neq i}f\left(
P_{k}\right) +f\left( P_{j}^{\prime}\right) \right) ,
\end{equation*}
where $\Delta_{i}\left( P_{j}^{\prime}\right) $ denotes the simplex obtained
from $\Delta$ by replacing the vertex $P_{i}$ by $P_{j}^{\prime}$. Notice
that (\ref{lambda}) yields $\lambda_{i}\left( P_{j}^{\prime}\right) \limfunc{%
Vol}\left( \Delta\right) =\limfunc{Vol}(\Delta _{i}\left(
P_{j}^{\prime}\right) ).$ By multiplying both sides by $\lambda_{i}\left(
P_{j}^{\prime}\right) $ and summing over $i$ we obtain%
\begin{align}
\frac{1}{\limfunc{Vol}\left( \Delta\right) }\int_{\Delta}f\left( x\right) 
\mathrm{d}x & \leq\frac{1}{n+1}\left( \sum_{i}\lambda_{i}\left(
P_{j}^{\prime}\right) \sum_{k\neq i}f\left( P_{k}\right) +f\left(
P_{j}^{\prime}\right) \right)  \notag \\
& =\frac{1}{n+1}\left( \sum_{i}\lambda_{i}\left( P_{j}^{\prime}\right)
\left( \sum_{k}f\left( P_{k}\right) -f\left( P_{i}\right) \right) +f\left(
P_{j}^{\prime}\right) \right)  \label{6} \\
& =\frac{1}{n+1}\left( \sum_{i}\left( 1-\lambda_{i}\left( P_{j}^{\prime
}\right) \right) f\left( P_{i}\right) +f\left( P_{j}^{\prime}\right) \right)
.  \notag
\end{align}
Furthermore, since $\Delta^{\prime}$ and $\Delta$ have the same barycenter,
we have 
\begin{equation*}
1-\lambda_{i}\left( P_{j}^{\prime}\right) =\sum_{k\neq j}\lambda_{i}\left(
P_{k}^{\prime}\right)
\end{equation*}
for all $i=1,...,n+1.$ Then 
\begin{equation*}
\frac{1}{\limfunc{Vol}\left( \Delta\right) }\int_{\Delta}f\left( x\right) 
\mathrm{d}x\leq\frac{1}{n+1}\left( \sum_{k\neq j}\sum_{i}\lambda_{i}\left(
P_{k}^{\prime}\right) f\left( P_{i}\right) +f\left( P_{j}^{\prime}\right)
\right) .
\end{equation*}
The fact that $\Delta^{\prime}$ and $\Delta$ have the same barycenter and
the convexity of the function $f$ yield 
\begin{multline*}
\sum_{j}f\left( P_{j}\right) =\sum_{k}\sum_{i}\lambda_{i}\left(
P_{k}^{\prime}\right) f\left( P_{i}\right) \\
=\sum_{k\neq j}\sum_{i}\lambda_{i}\left( P_{k}^{\prime}\right) f\left(
P_{i}\right) +\sum_{i}\lambda_{i}\left( P_{j}^{\prime}\right) f\left(
P_{i}\right) \\
\geq\sum_{k\neq j}\sum_{i}\lambda_{i}\left( P_{k}^{\prime}\right) f\left(
P_{i}\right) +f\left( \sum_{i}\lambda_{i}\left( P_{j}^{\prime}\right)
P_{i}\right) \, \\
=\sum_{k\neq j}\sum_{i}\lambda_{i}\left( P_{k}^{\prime}\right) f\left(
P_{i}\,\right) +f\left( P_{j}^{\prime}\,\right) .
\end{multline*}
This concludes the proof (\ref{R}).\newline
}

{\normalsize When $\Delta^{\prime}$ as a singleton, the inequality (\ref{R})
coincides with the inequality (\ref{15}). Notice that if we omit the
barycenter condition then the inequality (\ref{6}) translates into (\ref%
{rhh1}).\newline
}

{\normalsize An immediate consequence of Theorem \ref{thm3} is the following
result due to Farissi \cite{elf10}: }

{\normalsize COROLLARY 2. \label{cor2}Assume that $f$ : $\left[ a,b\right] $ 
$\rightarrow$ $\mathbb{R}$ is a convex function and $\lambda\in\left[ 0,1%
\right] $. Then%
\begin{align*}
\frac{1}{b-a}\int_{a}^{b}\,f(x)\,\mathrm{d}x & \leq\frac{\left(
1-\lambda\right) f\left( a\right) +\lambda f\left( b\right) +f\left( \lambda
a+\left( 1-\lambda\right) b\right) }{2} \\
& \leq\frac{f\left( a\right) +f\left( b\right) }{2}
\end{align*}
and%
\begin{align*}
\frac{1}{b-a}\int_{a}^{b}\,f(x)\,\mathrm{d}x & \geq\lambda f\left( \frac{%
a+\left( 1-\lambda\right) a+\lambda b}{2}\right) +\left( 1-\lambda\right)
f\left( \frac{b+\left( 1-\lambda\right) a+\lambda b}{2}\right) \\
& \geq f\left( \frac{a+b}{2}\right) .
\end{align*}
}

{\normalsize \textit{Proof}. Apply Theorem \ref{thm3} for $n=1$ and $%
\Delta^{\prime}$ the subinterval of endpoints $\left( 1-\lambda\right)
a+\lambda b\ $and $\lambda a+\left( 1-\lambda\right) b.$\newline
}

{\normalsize Another refinement of the Hermite-Hadamard inequality in the
case of simplices is as follows.\newline
}

{\normalsize THEOREM 4. \label{thm4}Let $\Delta\subset\mathbb{R}^{n}$ be an $%
n$-dimensional simplex of vertices $P_{1},...,P_{n+1}$ endowed with the
Lebesgue measure and let $\Delta^{\prime}\subseteq\Delta$ be a subsimplex
whose barycenter with respect to the normalized Lebesgue measure $\frac{%
\mathrm{d}x}{\limfunc{Vol}\left( \Delta ^{\prime}\right) }$ is $%
P=\sum_{k}\lambda_{k}\left( P\right) P_{k}.$ Then for every continuous
convex function $f:\Delta\rightarrow\mathbb{R},$%
\begin{equation}
f\left( P\right) \leq\frac{1}{\limfunc{Vol}\left( \Delta^{\prime }\right) }%
\int_{\Delta^{\prime}}\,f(x)\,\mathrm{d}x\leq\sum_{j}\lambda _{j}\left(
P\right) f\left( P_{j}\right) .  \label{12}
\end{equation}
}

{\normalsize \textit{Proof}. Let $P_{k}^{\prime},$ $k=1,...,n+1$ be the
vertices of $\Delta^{\prime}$. By Corollary \ref{cor1}, 
\begin{equation*}
f\left( P\right) \leq\frac{1}{\limfunc{Vol}\left( \Delta^{\prime }\right) }%
\int_{\Delta^{\prime}}\,f(x)\,\mathrm{d}x\leq\frac{1}{n+1}\sum _{k}f\left(
P_{k}^{\prime}\right) ,
\end{equation*}
The barycentric representation of each of the points $P_{k}^{\prime}\in%
\Delta $ gives us 
\begin{equation*}
\left\{ 
\begin{array}{l}
\sum_{j}\lambda_{j}\left( P_{k}^{\prime}\right) P_{j}\,=P_{k}^{\prime} \\ 
\sum_{k}\lambda_{k}\left( P_{j}^{\prime}\right) =1%
\end{array}
\right. .
\end{equation*}
Since $f$ is a convex function,%
\begin{align*}
\frac{1}{n+1}\sum_{k}f\left( P_{k}^{\prime}\right) & =\frac{1}{n+1}%
\sum_{k}f\left( \sum_{j}\lambda_{j}\left( P_{k}^{\prime}\right)
P_{j}\,\right) \\
& \leq\sum_{j}\left( \frac{1}{n+1}\sum_{k}\lambda_{j}\left( P_{k}^{\prime
}\right) \right) f\left( P_{j}\,\right) =\sum_{j}\lambda_{j}\left( P\right)
f\left( P_{j}\right)
\end{align*}
and the assertion of Theorem \ref{thm4} is now clear.\newline
}

{\normalsize As a corollary of Theorem \ref{thm4} we get the following
result due Vasi\'{c} and Lackovi\'{c} \cite{vas76}, and Lupa\c{s} \cite%
{lup76} (cf. J. E. Pe\v{c}ari\'{c} et al. \cite{pec92}).\newline
}

{\normalsize COROLLARY 3. \label{cor3}Let $p$ and $q$ be two positive
numbers and $a_{1}\leq a\leq b\leq b_{1}.$ Then the inequalities%
\begin{equation*}
f\left( \frac{pa+qb}{p+q}\right) \leq\frac{1}{2y}\int_{A-y}^{A+y}\,f(x)\,%
\mathrm{d}x\leq\frac{pf\left( a\right) +qf\left( b\right) }{p+q}
\end{equation*}
hold for $A=\frac{pa+qb}{p+q},$ $y>0$ and all continuous convex functions $f:%
\left[ a_{1},b_{1}\right] \rightarrow\mathbb{R}$ if and only if 
\begin{equation*}
y\leq\frac{b-a}{p+q}\min\left\{ p,q\right\} .
\end{equation*}
}

{\normalsize The right-hand side of the inequality stated in Theorem \ref%
{thm4} can be improved as follows.\newline
}

{\normalsize THEOREM 5. \label{Thm5}Suppose that $\Delta\ $is an $n$%
-dimensional simplex of vertices $P_{1},...,P_{n+1}$ and let $P\in\Delta$.
Then for every subsimplex $\Delta^{\prime}\subset\Delta$ such that $%
P=\sum_{j}\lambda_{j}\left( P\right) P_{j}=b_{\mathrm{d}x/\limfunc{Vol}%
\left( \Delta^{\prime }\right) }$ we have%
\begin{equation*}
\frac{1}{\limfunc{Vol}\left( \Delta^{\prime}\right) }\int
_{\Delta^{\prime}}\,f(x)\,\mathrm{d}x\leq\frac{1}{n+1}\left(
n\sum_{j}\lambda_{j}\left( P\right) f\left( P_{j}\right) +f\left( P\right)
\right)
\end{equation*}
for every continuous convex function $f:\Delta\rightarrow\mathbb{R}\mathbf{.}
$ }

{\normalsize \textit{Proof}. Let $P_{k}^{\prime},$ $k=1,...,n+1$ be the
vertices of $\Delta^{\prime}$. We denote by $\Delta_{i}^{\prime}$ the
subsimplex obtained by replacing the vertex $P_{i}^{\prime}$ by the
barycenter $P$ of the normalized measure $\mathrm{d}x/\limfunc{Vol}\left(
\Delta^{\prime}\right) $ on $\Delta^{\prime}$. }

{\normalsize According to Corollary \ref{cor1}, 
\begin{align*}
\frac{1}{\limfunc{Vol}\left( \Delta_{i}^{\prime}\right) }\int
_{\Delta_{i}^{\prime}}\,f(x)\,\mathrm{d}x & \leq\frac{1}{n+1}\sum_{k\neq
i}f\left( P_{k}^{\prime}\right) +\frac{1}{n+1}f\left( P\right) \\
& =\frac{1}{n+1}\sum_{k\neq i}f\left( \sum_{j}\lambda_{j}\left(
P_{k}^{\prime}\right) P_{j}\,\right) +\frac{1}{n+1}f\left( P\right) \\
& \leq\sum_{j}\left( \frac{1}{n+1}\sum_{k\neq i}\lambda_{j}\left(
P_{k}^{\prime}\right) \right) f\left( P_{j}\,\right) +\frac{1}{n+1}f\left(
P\right) \\
& =\sum_{j}\left( \lambda_{j}\left( P\right) -\frac{1}{n+1}\lambda
_{j}\left( P_{i}^{\prime}\right) \right) f\left( P_{j}\,\right) +\frac {1}{%
n+1}f\left( P\right)
\end{align*}
for each index $i=1,...,n+1.$ Summing up over $i$ we obtain%
\begin{align*}
& \frac{n+1}{\limfunc{Vol}\left( \Delta^{\prime}\right) }\int
_{\Delta^{\prime}}\,f(x)\,\mathrm{d}x \\
& \leq\left( n+1\right) \sum_{j}\lambda_{j}\left( P\right) f\left(
P_{j}\,\right) -\sum_{i}\frac{1}{n+1}\sum_{j}\lambda_{j}\left(
P_{i}^{\prime}\right) f\left( P_{j}\,\right) +f\left( P\right) \\
& =n\sum_{j}\lambda_{j}\left( P\right) f\left( P_{j}\,\right) +f\left(
P\right) .
\end{align*}
and the proof of the theorem is completed.\newline
}

{\normalsize Of course, the last theorem yields an improvement of Corollary %
\ref{cor3}. This was first noticed in \cite[pp. 146]{pec92}. }

{\normalsize We end this paper with an alternative proof of some particular
case of a result established by A. Guessab and G. Schmeisser [4, Theorem
2.4]. Before we start we would like to turn the reader's attention to the
paper [12] by S. W\k asowicz, where the result we are going to present was
obtained by some more general approach (see [11, Theorems 2 and Corollary
4]). \newline
}

{\normalsize THEOREM 6. \label{thm6}Let $\Delta\subset\mathbb{R}^{n}$ be an $%
n$-dimensional simplex of vertices $P_{1},...,P_{n+1}$ endowed with the
Lebesgue measure and let $M_{1},...,M_{m}$ be points in $\Delta$ such that $%
b_{\mathrm{d}x/\limfunc{Vol}\left( \Delta\right) }$ is a convex combination $%
\sum_{j}\beta_{j}M_{j}$ of them. Then for every continuous convex function $%
f:\Delta\rightarrow\mathbb{R},$ the following inequalities hold:%
\begin{equation*}
f\left( b_{\mathrm{d}x/\limfunc{Vol}\left( \Delta\right) }\right)
\leq\sum_{j=1}^{m}\beta_{j}f\left( M_{j}\right) \leq\frac{1}{n+1}\sum
_{i=1}^{n+1}f\left( P_{i}\right) .
\end{equation*}
}

{\normalsize \textit{Proof}. The first inequality follows from Jensen's
inequality. In order to establish the second inequality, since we have 
\begin{equation*}
b_{\mathrm{d}x/\limfunc{Vol}\left( \Delta\right) }=\sum_{j}\beta _{j}\left(
\sum_{i}\lambda_{i}\left( M_{j}\right) P_{i}\right) =\sum _{i}\left(
\sum_{j}\beta_{j}\lambda_{i}\left( M_{j}\right) \right) P_{i}=\frac{1}{n+1}%
\sum_{i}P_{i}\,,
\end{equation*}
we infer that $\sum_{j}\beta_{j}\lambda_{i}\left( M_{j}\right) =\frac {1}{n+1%
}$ for every $i$ and 
\begin{align*}
\sum_{j}\beta_{j}f\left( M_{j}\right) & =\sum_{j}\beta_{j}f\left(
\sum_{i}\lambda_{i}\left( M_{j}\right) P_{i}\,\right) \leq\sum_{i}\left(
\sum_{j}\beta_{j}\lambda_{i}\left( M_{j}\right) \right) f\left(
P_{i}\,\right) \\
& =\frac{1}{n+1}\sum_{i}f\left( P_{i}\,\right) .
\end{align*}
This completes the proof.\newline
}

{\normalsize \textbf{Acknowledgement.} The first author has been supported
by CNCSIS Grant $420/2008$. The authors are indebted to Szymon W\k asowicz
for pointing out to them the references \cite{was10-b} and \cite{gue}. }

{\normalsize \providecommand{\bysame}{\leavevmode\hbox
to3em{\hrulefill}\thinspace} }

\begin{flushright}
{\normalsize {\scriptsize \textit{Department of Mathematics\newline
University of Craiova\newline
Street A. I. Cuza 13, Craiova, RO-200585 \newline
Romania\newline
fcmitroi@yahoo.com, catalin\_gsep@yahoo.com} } }
\end{flushright}

\end{document}